\documentclass{elsart3-1}

\usepackage{lmodern}  
\usepackage[T1]{fontenc}
\usepackage{amsmath}
\usepackage{amssymb}
\providecommand{\qedsymbol}{$\square$}
\newcommand{\mathqed}{\quad\hbox{\qedsymbol}}
\DeclareRobustCommand{\qed}{%
  \ifmmode \mathqed
  \else
    \leavevmode\unskip\penalty9999 \hbox{}\nobreak\hfill
    \quad\hbox{\qedsymbol}%
  \fi
}
\usepackage{graphicx}
\usepackage{mathrsfs}
\usepackage[english,francais]{babel}

\newcommand{\C}{\mathbb{C}}

\newcommand{\A}{\mathbb{A}}

\newcommand{\lf}{\mathbb{L}}

\newcommand{\rk}{\textnormal{rk}}
\newtheorem{theorem}{Theorem}[section]
\newtheorem{lemma}[theorem]{Lemma}
\newtheorem{e-proposition}[theorem]{Proposition}

\newtheorem{e-definition}[theorem]{Definition\rm}

\renewcommand{\theequation}{\arabic{equation}}
\def\og{\leavevmode\raise.3ex\hbox{$\scriptscriptstyle\langle\!\langle$~}}
\def\fg{\leavevmode\raise.3ex\hbox{~$\!\scriptscriptstyle\,\rangle\!\rangle$}}

\journal{HAL}
\begin{document}
\centerline{}
\begin{frontmatter}


\selectlanguage{english}
\title{The class of the affine line is a zero divisor in the Grothendieck ring: an improvement}


\vspace{0.8cm}

\selectlanguage{english}
\author{Nicolas Martin},
\ead{nicolas.martin@polytechnique.edu}
\address{Centre de Math\'ematiques Laurent Schwartz, École polytechnique, 91128 Palaiseau cedex, France}
\date{April 17, 2016}


\medskip

\vspace{1cm}

\begin{abstract}
\selectlanguage{english}
\indent
Lev A. Borisov has shown that the class of the affine line is a zero divisor in the Grothendieck ring of algebraic varieties over complex numbers. We improve the final formula by removing a factor.

\vskip 2\baselineskip

\selectlanguage{francais}
\noindent{\bf R\'esum\'e}
\vskip 0.5\baselineskip
\setlength\parindent{12pt}
Lev A. Borisov a prouv\'e que la classe de la droite affine est un diviseur de z\'ero dans l'anneau de Grothendieck des vari\'et\'es alg\'ebriques complexes. Nous am\'eliorons la formule finale en supprimant un facteur.

\end{abstract}
\end{frontmatter}

\selectlanguage{english}

\vspace{0.5cm}

\section{Introduction}

The Grothendieck ring $K_0(\textnormal{Var}_{\C})$ of complex algebraic varieties is defined as the quotient of the free abelian group generated by the isomorphism classes $[X]$ of complex algebraic varieties modulo the relations
$$[X]=[Y]+[X \setminus Y]$$
for all closed subvarieties $Y \subset X$. The cartesian product of varieties gives the product structure.

\newpage

The class $\lf=[\A^1(\C)]$ of the affine line has a major role in the study of the Grothendieck ring. It has been proved in \cite{LL03} that $X$ and $Y$ are stably birational if and only if their classes $[X]$ and $[Y]$ are equal modulo $\lf$. After Bjorn Poonen had shown in \cite{Poo02} that $K_0(\textnormal{Var}_{\C})$ is not a domain, Lev Borisov has made precise this result in \cite{Bor14} by showing that $\lf$ is a zero divisor. He has compared the two sides $[X_W]$ and $[Y_W]$ of the Pfaffian-Grassmannian double mirror correspondence, and obtained the following formula:
$$([X_W]-[Y_W]) \cdot (\lf^2-1) \cdot (\lf-1) \cdot \lf^7=0.$$

\vspace{0.5cm}

This result is not only an improvement of that of Poonen: it is crucial in motivic integration to understand the kernel of the localization morphism $K_0(\textnormal{Var}_{\C}) \rightarrow K_0(\textnormal{Var}_{\C})[\lf^{-1}]$, since we consider classes in the localized ring. In this paper, we improve this formula as follows.

\vspace{0.5cm}

\begin{theorem}
$([X_W]-[Y_W]) \cdot \lf^6=0$
\end{theorem}

\vspace{0.5cm}

\noindent
\textbf{Acknowledgements.} The author is indebted to Johannes Nicaise and Claude Sabbah for their careful reading and constructive comments on the preliminary version of the note written in September of 2015. The note has particularly benefited from relevant comments of Antoine Chambert-Loir in November of 2015.

\section{The class of Grasmannians}

\begin{e-proposition}
For $2 \leq k < n$, we have the relation
$$[G(k,n)]=[G(k,n-1)]+\lf^{n-k} \cdot [G(k-1,n-1)].$$
\end{e-proposition}

\vspace{0.4cm}

\noindent
\textbf{Proof.} Let $e_1,...,e_n$ be the canonical basis of $\C^n$, $F$ the hyperplane orthogonal to $e_n$, $U \subset G(k,n)$ the open subset defined by $\{ T \in G(k,n) \ | \ \dim(T \cap F)=k-1 \}$ and $\pi : U \rightarrow G(k-1,F)$ the regular mapping which sends $T$ on $T \cap F$. For $S \in G(k-1,F)$, the fiber $\pi^{-1}(S)$ can be identified to
$$\mathbb{P}(\C^n/S) \setminus \mathbb{P}(F/S) \simeq \A^{n-k}.$$
Let $H$ be a complementary subspace of $S$ in $F$ and the open subset $V=\{ S' \in G(k-1,F) \ | \ S' \oplus H = F \}$. For all $S' \in V$, we have the identification $\C^n/S' \simeq H \oplus \C e_n$, hence $\pi$ is a trivial fibration over $V$. Consequently, $\pi$ is a locally trivial fibration, therefore $[U]=\lf^{n-k} \cdot [G(k-1,n-1)]$. We have $[G(k,n)]=[Z]+[U]$ with $Z=G(k,n) \setminus U=\{ T \in G(k,n) \ | \ T \subset F \}=G(k,F)$, which shows the announced formula. \qed

\vspace{0.5cm}

A simple induction gives the following formulas for $n \geq 4$:
$$[G(2,n)]=\left\{
\begin{array}{ll}
  \left[\mathbb{P}^{n-2}\right] \cdot \sum\limits_{k=0}^{(n-2)/2} \lf^{2k} & \textnormal{if } n \textnormal{ is even}\\[10pt]
  \left[\mathbb{P}^{n-1}\right] \cdot \sum\limits_{k=0}^{(n-3)/2} \lf^{2k} & \textnormal{if } n \textnormal{ is odd.}
\end{array}
\right.$$

For example, $[G(2,5)]=[\mathbb{P}^4] \cdot (\lf^2+1)$ and $[G(2,7)]=[\mathbb{P}^6] \cdot (\lf^4+\lf^2+1)$.

\newpage

\section{Improvement of Borisov's formula}

\subsection{Pfaffian and Grassmannian double mirror varieties}

Let $V$ be a 7-dimensional complex vector space and $W$ a generic 7-dimensional space of skew forms on $V$. We define $X_W$ as a subvariety of the Grassmannian $G(2,V)$ which is the locus of all $T \in G(2,V)$ with $\omega_{|T}=0$ for all $\omega \in W$, and $Y_W$ as a subvariety of $\mathbb{P}W$ of skew forms whose rank is less than 6. Smoothness of these two varieties has been shown by E. Rødland in \cite{Rod00}. Furthermore, we know that all forms in $Y_W$ have rank 4 and all forms in $\mathbb{P}W \setminus Y_W$ have rank 6.

\subsection{The formula}

Let us define $H$ as a subvariety of $G(2,V) \times \mathbb{P}W$ which consists of pairs $(T,\C \omega)$ with $\omega_{|T}=0$. In order to obtain the explicit equations which define $H$, let us set $T_0 \in G(2,V)$ with basis $e_1,e_2$ and $H$ a complementary subspace with basis $e_3,...,e_7$. The neighborhood $U=\{T \in G(2,V) \ | \ T \oplus H=V\}$ of $T_0$ can be identified to $\mathscr{L}(T_0,H)$ by considering the map $f \in \mathscr{L}(T_0,H) \mapsto \{x+f(x) \ | \ x \in T_0\} \in U$. If we set $(f_{i,j})_{(i,j)\in \{1,2\}\times\{3,...,7\}}$ the basis of $\mathscr{L}(T_0,H)$ adapted to the two bases previously considered, we can identify $T \in U$ to $\{ x + \sum\alpha_{i,j}f_{i,j}(x) \ | \ x \in T_0\}$. Now, for $\omega=\sum_{i=1}^{7}\beta_i \omega_i \in W$, the condition $\omega_{|T}=0$ can be expressed as
$$\sum_{i=1}^{7}\beta_i \omega_i\left(e_1+\sum_{j=3}^{7}\alpha_{1,j}e_j,e_2+\sum_{j=3}^{7}\alpha_{2,j}e_j\right)=0.$$

\vspace{0.3cm}

Looking at the projections onto the two factors $G(2,V)$ and $\mathbb{P}W$ will give us two ways to express $[H]$. Theorem 1.1 will be a direct consequence of the two next propositions.

\vspace{0.5cm}

\begin{e-proposition}
$[H]=[\mathbb{P}^6] \cdot (\lf^4+\lf^2+1) \cdot [\mathbb{P}^5]+[X_W] \cdot \lf^6$
\end{e-proposition}

\vspace{0.5cm}

\noindent
\textbf{Proof.} Considering the projection $p : H \rightarrow G(2,V)$ onto the first factor, which is a trivial fibration in restriction to $p^{-1}(X_W)$ and a locally trivial fibration in restriction to $G(2,V) \setminus p^{-1}(X_W)$, Proposition 2.4 of \cite{Bor14} proves that
$$[H]=[G(2,7)] \cdot [\mathbb{P}^5]+[X_W] \cdot \lf^6.$$

\vspace{0.2cm}

\noindent
The expression $[G(2,7)]=[\mathbb{P}^6] \cdot (\lf^4+\lf^2+1)$ gives the result. \qed

\vspace{0.5cm}

\begin{e-proposition}
$[H]=[Y_W] \cdot \lf^6 + [\mathbb{P}^6] \cdot [\mathbb{P}^5] \cdot (\lf^4+\lf^2+1)$
\end{e-proposition}

\vspace{0.5cm}

\begin{lemma}
Let $\pi : H \rightarrow \mathbb{P}W$ be the projection onto the second factor. Its restrictions to $\pi^{-1}(Y_W)$ and $\pi^{-1}(\mathbb{P}W \setminus Y_W)$ are piecewise trivial fibrations (see 4.2.1 in \cite{Seb04}).
\end{lemma}

\vspace{0.3cm}

\noindent
\textbf{Proof of the lemma.} The reasoning is the same for rank 4 ($Y_4=Y_W$) and rank 6 ($Y_6=\mathbb{P}W \setminus Y_W$). For $i \in \{4,6\}$, let us set 
$$Z_i=\pi^{-1}(Y_i)=H \cap (G(2,V) \times Y_i).$$

\vspace{0.2cm}

\noindent
In order to have piecewise triviality of $\pi$ on $Z_i$, it suffices, according to Theorem 4.2.3 in \cite{Seb04}, to prove that there exists a uniform fiber $F_i$ such that for all $x \in Y_i$,
$$Z_i \times_{Y_i} \{x\} \simeq F_i \times_\C \textnormal{Spec}(\kappa(x)).$$
\noindent
To achieve this, it suffices to note that a skew form of rank 4 or 6 with coefficients in a field $K \supset \C$ is congruent to the skew form
$$\begin{pmatrix}
   \phantom{-}0 & \phantom{-}I_2 & \phantom{-}0 \phantom{-}\\
   -I_2 & \phantom{-}0 & \phantom{-}0 \phantom{-}\\
   \phantom{-}0 & \phantom{-}0 & \phantom{-}0 \phantom{-}
\end{pmatrix} \quad \textnormal{or} \quad
\begin{pmatrix}
   \phantom{-}0 & \phantom{-}I_3 & \phantom{-}0 \phantom{-}\\
   -I_3 & \phantom{-}0 & \phantom{-}0 \phantom{-}\\
   \phantom{-}0 & \phantom{-}0 & \phantom{-}0 \phantom{-}
\end{pmatrix}$$

\vspace{0.2cm}
\noindent
with a base change having coefficients in $K$, an action that spreads on fibers. \qed

\vspace{0.5cm}

\begin{lemma}
Let $\C \omega \in Y_W$ be a closed point. Then the class of its fiber is
$$[\pi^{-1}(\C \omega)]=[\mathbb{P}^5] \cdot (\lf^4+\lf^2+1)+\lf^6.$$
\end{lemma}

\vspace{0.3cm}

\noindent
\textbf{Proof.} As $\rk(\omega)=4$, there exists a basis $e_1,...,e_7$ of $V$ in which the matrix of $\omega$ is
$$\begin{pmatrix}
   \phantom{-}0 & \phantom{-}I_2 & \phantom{-}0 \phantom{-}\\
   -I_2 & \phantom{-}0 & \phantom{-}0 \phantom{-}\\
   \phantom{-}0 & \phantom{-}0 & \phantom{-}0 \phantom{-}
\end{pmatrix}.$$

\noindent
Denote $F=\textnormal{Vect}\{e_3,...,e_7\}$ and $H=F \oplus \C e_2$. We have
$$[\pi^{-1}(\C \omega)]=[\{T \in G(2,V) \ | \ \omega_{|T}=0\}]=[\{T \in G(2,H) \ | \ \omega_{|T}=0\}]+[U]$$
where $U$ is the open subset $\{T \in G(2,V) \ | \ \dim(T \cap H)=1, \ \omega_{|T}=0\}$, with the locally trivial fibration $\pi : U \rightarrow \mathbb{P}H=\mathbb{P}^5$. Note that $\ker(\omega)=\textnormal{Vect}\{e_5,e_6,e_7\} \subset H$ and $\ker(\omega_{|H})=\ker(\omega) \oplus \C e_3 \subset H$.

\vspace{0.3cm}

\noindent
Let $D = \C e \in \mathbb{P}H$. There are three cases.

\vspace{0.3cm}

\noindent
$\bullet$ First case: $D \subset \ker(\omega)$. We have
\begin{align*}
[\pi^{-1}(D)] &= [\{ \C f \in \mathbb{P}(V/D) \ | \ \omega(f,e)=0 \}]-[\{ \C f \in \mathbb{P}(H/D) \ | \ \omega_{|H}(f,e)=0 \}]\\
&= [\mathbb{P}^{5}]-[\mathbb{P}^4]=\lf^5.
\end{align*}

\vspace{0.3cm}

\noindent
$\bullet$ Second case: $D \not\subset \ker(\omega)$ and $D \subset \ker(\omega_{|H})$. In this case $\pi^{-1}(D)=\varnothing$, because
$$\{ \C f \in \mathbb{P}(V/D) \ | \ \omega(f,e)=0 \}=\{ \C f \in \mathbb{P}(H/D) \ | \ \omega_{|H}(f,e)=0 \}.$$

\vspace{0.3cm}

\noindent
$\bullet$ Third case: $D \not\subset \ker(\omega_{|H})$. We have
\begin{align*}
[\pi^{-1}(D)] &= [\{ \C f \in \mathbb{P}(V/D) \ | \ \omega(f,e)=0 \}]-[\{ \C f \in \mathbb{P}(H/D) \ | \ \omega_{|H}(f,e)=0 \}]\\
&= [\mathbb{P}^{4}]-[\mathbb{P}^3]=\lf^4.
\end{align*}

\vspace{0.2cm}

\noindent
Consequently
\begin{align*}
[U] &= [\mathbb{P}\ker(\omega)] \cdot \lf^5+([\mathbb{P}H]-[\mathbb{P}\ker(\omega_{|H})]) \cdot \lf^4\\
&= [\mathbb{P}^2] \cdot \lf^5+([\mathbb{P}^5]-[\mathbb{P}^3]) \cdot \lf^4\\
&= ([\mathbb{P}^5]-1) \cdot \lf^4.
\end{align*}

\vspace{0.2cm}

\noindent
We can repeat the argument with $H$. As $\omega_{|F}=0$, we have
\begin{align*}
[\{T \in G(2,H) \ | \ \omega_{|T}=0\}] &= [\{T \in G(2,F) \ | \ \omega_{|T}=0\}]+[\mathbb{P}\ker(\omega_{|H})] \cdot \lf^4\\
&= [G(2,5)]+[\mathbb{P}^3] \cdot \lf^4\\
&= [\mathbb{P}^4] \cdot (\lf^2+1)+[\mathbb{P}^3] \cdot \lf^4.
\end{align*}

\noindent
Finally, we get
\begin{align*}
[\pi^{-1}(\C \omega)] &= ([\mathbb{P}^5]-1) \cdot \lf^4+[\mathbb{P}^4] \cdot (\lf^2+1)+[\mathbb{P}^3] \cdot \lf^4\\
&= ([\mathbb{P}^5]-1) \cdot \lf^4+([\mathbb{P}^5]-\lf^5) \cdot (\lf^2+1)+ (\lf^3 + \lf^2 + \lf + 1) \cdot \lf^4\\
&= [\mathbb{P}^5] \cdot (\lf^4+\lf^2+1)+\lf^6.
\end{align*}\qed

\vspace{0.5cm}

A similar calculation gives the following result.

\vspace{0.5cm}

\begin{lemma}
Let $\C \omega \in \mathbb{P}W \setminus Y_W$ be a closed point. Then the class of its fiber is
$$[\pi^{-1}(\C \omega)] = [\mathbb{P}^5] \cdot (\lf^4+\lf^2+1).$$
\end{lemma}

\vspace{0.5cm}

\noindent
\textbf{Proof of Proposition 3.2.} Let $\C \omega_1 \in Y_W$ and $\C \omega_2 \in \mathbb{P}W \setminus Y_W$ be two closed points. Lemma 3.3 implies that

$$\left\{
\begin{array}{ll}
[\pi^{-1}(Y_W)]=[Y_W] \cdot [\pi^{-1}(\C \omega_1)]\\[5pt]
[\pi^{-1}(\mathbb{P}W \setminus Y_W)]=([\mathbb{P}W]-[Y_W]) \cdot [\pi^{-1}(\C \omega_2)],
\end{array}
\right.$$

\noindent
and consequently

$$[H] = [Y_W] \cdot [\pi^{-1}(\C \omega_1)]+([\mathbb{P}W]-[Y_W]) \cdot [\pi^{-1}(\C \omega_2)].$$

\vspace{0.5cm}

\noindent
Using Lemmas 3.4 and 3.5, we have
\begin{align*}
[H]&= [Y_W] \cdot ([\mathbb{P}^5] \cdot (\lf^4+\lf^2+1)+\lf^6)+([\mathbb{P}^6]-[Y_W]) \cdot [\mathbb{P}^5] \cdot (\lf^4+\lf^2+1)\\
&= [Y_W] \cdot \lf^6 + [\mathbb{P}^6] \cdot [\mathbb{P}^5] \cdot (\lf^4+\lf^2+1),
\end{align*}

\noindent
which concludes the proof. \qed

\vspace{1.2cm}

\begin{center}
\begin{small}
REFERENCES
\end{small}
\end{center}

\vspace{-1cm}

\renewcommand{\refname}{}
\bibliographystyle{amsalpha}
\bibliography{bibli}

\end{document}